\documentclass{am2013}

\usepackage{epsfig}
\usepackage{float}
\usepackage{amsmath}
\usepackage{setspace} 
\makeatother
\usepackage{makecell}
\usepackage{longtable}
\usepackage{tabularx}
\usepackage{changepage}
\newlength{\extralength}
\usepackage{array}
\usepackage{tabularx}
\usepackage{algorithm}
\usepackage{algorithmic}
\newcolumntype{C}{>{\centering\arraybackslash}X}
\usepackage{booktabs}
\usepackage{caption}
\title{RBF Kernel Parameter Formula for Data Classification Methods}
\newlength{\fulllength}
\author{Lakhdar Remaki\\
\vskip 2mm {\small
Departement of Mathematics and Computer Science, \\
Alfaisal University, KSA\\
lremaki@alfaisal.edu
}
}

\begin{document}
 
\maketitle

\begin{abstract}
Radial Basis Function (RBF), or Gaussian, kernels are among the most widely used parametric kernels in machine learning, particularly in methods such as Support Vector Machines (SVM) and kernel-based subspace approaches. The kernel parameter $\gamma$ (or $\sigma$ in the Gaussian formulation) must be carefully tuned, as the performance of these methods strongly depends on its value and is highly sensitive to improper selection. In practice, this parameter is typically determined through computationally expensive training procedures, which may also lack robustness.
In this paper, we propose an efficient analytical formula for selecting the RBF kernel parameter that significantly reduces the computational cost of RBF-based methods. The proposed approach is derived by optimizing the diameter of mapped classes in the feature space while simultaneously maximizing inter-class feature distances. The detailed formulation is presented, and its efficiency is validated on the widely used SVM algorithm as well as on a Proper Orthogonal Decomposition (POD)-based subspace method for both binary and multi-class classification problems.
\end{abstract}

\section*{keyword}
RBF, $\gamma$ parameter, SVM, KOS, kernels


\section{Introduction}
 
Machine Learning (ML) is a prominent and rapidly evolving research area at the core of Artificial Intelligence (AI), with a wide range of applications, including handwriting recognition, automated disease detection, robotics, and others. Kernel methods constitute an important class of data classification techniques that have demonstrated high efficiency across many applications. The fundamental idea of kernel methods is to map the original data, referred to as attributes, into a higher-dimensional space, known as the feature space, where the data can be more easily separated (i.e., become linearly separable).

In practice, the mapping function is not explicitly defined, as this is a tedious task; instead, a suitable kernel function is used, provided that all operations can be expressed in terms of dot products. The existence of such kernels is guaranteed under certain conditions by Mercer’s theorem~\cite{Merc909,CorVap95}. The performance of these methods depends strongly on the choice of the kernel. Various kernels have been proposed in the literature, and research on this topic remains highly active.

Radial Basis Function (RBF) kernels are among the most widely used parametric kernels in machine learning, particularly in Support Vector Machines (SVM) (see~\cite{Merc909,CorVap95,BWang2021,Borah020,Gaye021,Tian012,KMOD001,Wang024,Wang023,Shao023} for recent developments) and kernel-based subspace learning methods, such as Kernel Principal Component Analysis (KPCA)~\cite{kpca1997,ksubspace2022} and other subspace approaches~\cite{lremaki2026, subspace2016,subspace2022,subspace2023}.

One of the theoretical reason its effectiveness is that the RBF kernel induces a mapping of the original data into an infinite-dimensional feature space, thereby enhancing data separability. However, the performance of the RBF kernel critically depends on the choice of the kernel hyperparameter $\gamma$ (or equivalently $\sigma$ in the Gaussian formulation). An improper selection may severely degrade generalization performance. In practice, this parameter is typically determined by learning via cross-validation or grid search strategies~\cite{Shekar19, Lin2008, Syarif2016, Bergstra12, Hinton06}, which can lead to substantial computational overhead and limited robustness, especially for large-scale or multi-class problems.

In this work, we introduce a computationally efficient analytical formulation for the direct estimation of the RBF kernel parameter. The proposed formula is derived from geometric considerations in the induced feature space, where we jointly optimize intra-class compactness, by minimizing the mapped class diameter, and inter-class separability, by maximizing pairwise feature-space distances. This leads to a closed-form expression that eliminates the need for exhaustive hyperparameter search.

The derivation of the proposed formula is presented in detail, and its performance is demonstrated through multi-class test cases. The formula is evaluated using both SVM and Kernel Optimization Subspaces (KOS)~\cite{lremaki2026,lremaki2024}. Results obtained with the proposed method are compared to those from the classical learning process. The LIBSVM software package~\cite{LIBSVM} is used to evaluate the SVM method, while a homemade implementation is used for KOS. Experiments are conducted on selected real-world datasets from the LIBSVM~\cite{LIBSVM} and OpenML~\cite{openML} repositories.

The results show that the proposed formula achieves performance comparable to that of the standard learning process, while significantly reducing computational time.

    In \mbox{Section \ref{2}}, a~short overview of SVM and KOS is provided first, then details of the proposed $\gamma$ formula are developed in \mbox{Section \ref{3}}. Tests to validate the proposed formula are performed and reported in \mbox{Section \ref{4}}. Conclusions are drawn in \mbox{Section \ref{5}}.

\section{Short Overview of SVM and KOS Methods}\label{2}
\subsection{SVM Method}
Support Vector Machines (SVM) are supervised machine learning algorithms widely used for classification and regression tasks. The main objective of SVM is to determine an optimal decision boundary, called a \textit{hyperplane}, that separates data points belonging to different classes with the maximum possible margin. The margin is defined as the distance between the hyperplane and the nearest data samples from each class, known as \textit{support vectors}. Maximizing this margin improves the generalization capability of the model.

For a binary classification problem with a training set $(x_i, y_i)$, where $x_i$ represents the feature vector and $y_i \in \{-1,1\}$ denotes the class label, the hyperplane is defined as:

\begin{equation}
w \cdot x + b = 0
\end{equation}

where $w$ is the weight vector and $b$ is the bias term. The optimal hyperplane is obtained by solving the following optimization problem:

\begin{equation}
\min_{w,b} \frac{1}{2} ||w||^2
\end{equation}

subject to

\begin{equation}
y_i (w \cdot x_i + b) \geq 1, \quad \forall i
\end{equation}

When the data are not linearly separable, slack variables $\xi_i \geq 0$ are introduced to allow some classification errors, leading to the soft-margin formulation:

\begin{equation}
\min_{w,b,\xi} \frac{1}{2} ||w||^2 + C \sum_{i=1}^{n} \xi_i
\end{equation}

subject to

\begin{equation}
y_i (w \cdot x_i + b) \geq 1 - \xi_i, \quad \forall i
\end{equation}

To handle nonlinear classification problems, SVM employs the \textit{kernel trick}. Instead of working directly in the input space, the data are implicitly mapped into a higher-dimensional feature space via a nonlinear mapping $\phi(x)$. The key idea is that the optimization problem depends only on inner products between data points, which can be replaced by a kernel function:

\begin{equation}
K(x_i, x_j) = \langle \phi(x_i), \phi(x_j) \rangle
\end{equation}

This allows the algorithm to operate in the high-dimensional feature space without explicitly computing $\phi(x)$, significantly reducing computational complexity. Common kernel functions include the linear kernel, polynomial kernel, radial basis function (RBF) kernel, and sigmoid kernel. Through this approach, SVM can construct highly nonlinear decision boundaries in the original input space.

Due to its strong theoretical foundation and good generalization performance, SVM has been successfully applied in various domains such as image recognition, text classification, bioinformatics, and pattern recognition.
\subsection{KOS Method}
The KOS method is a kernel-based subspace approach in which feature subspaces are constructed using the Proper Orthogonal Decomposition (POD) technique \cite{volkwein2013proper}. An unseen sample is classified, as commonly done in standard subspace-based methods, by projecting its feature vector onto each class subspace and assigning the label corresponding to the minimum distance. Unlike SVM, KOS replaces separating hyperplanes, which can be a source of practical limitations, with optimal class subspaces that directly model the feature data. The method can be summarized as follows; for more details, see \cite{lremaki2026}.

\subsection{Summary of the~Algorithm}
The KOS method is summarized below, and its main computational steps are outlined.

Let $C_{k,~k=1,P}$ denote $P$ attribute sets representing $P$ different classes used during the learning stage. Let $B_{k=1,P}$ denote the corresponding mapped sets in the feature space. The~KOS classification procedure is as follows:
\begin{enumerate}
\item Check whether the classes $C_{k}$ are imbalanced. If~so, split the larger classes (roughly classes with size twice the size of the smallest class). For~classification, if~the minimum distance is obtained from a POD subspace of a subdivided class, the~unseen sample is assigned to that class. 
\item Select a \textit{Mercer} kernel $K$, then compute the kernel matrices $K_{k}$ for each feature class $C_{k}$ (of feature sub-classes if applicable). 
\begin{equation*}
K_{k}(i,j) = K(X_{i},Y_{j}).
\end{equation*}

\item Compute the POD correlation matrices $M_{k}$ for each feature class (or sub-class)  $C_{k}$ 
\begin{equation}\label{kernel-cent}
M_{k}(i,j)= K_{k}(i,j) -\frac{1}{N}\sum_{k=1}^{N}(K_{k}(i,k)+K_{k}(j,k)) + \frac{1}{N^{2}}\sum_{k=1}^{N}\sum_{l=1}^{N}K_{k}(l,k).
\end{equation}
\item Compute the eigenvalues and eigenvectors of each matrix $M_{k,~k=1,P}$.
\item For each unseen vector $\widehat{X}$, compute its coordinates $\alpha^{k}_{i}, i=1,size(C_{k})$ in each POD feature subspace $\phi(C_{k})$:

\begin{align}\label{coord-cent}
\alpha^{k}_{i}= \dfrac{1}{\sqrt{\sigma}_{i}}\bigg[\sum_{l=1}^{N}V^{i}_{l}K(X_{l},\widehat{X})  -\dfrac{1}{N} \sum_{l=1}^{N}\sum_{j=1}^{N}V^{i}_{l}K_{l,j} \\ \nonumber - (\frac{1}{N}\sum_{j=1}^{N}K(X_{j},\widehat{X}))\sum_{l=1}^{N}V^{i}_{l} +(\frac{1}{N^{2}} \sum_{k=1}^{N}\sum_{j=1}^{N}K_{l,j})\sum_{k=1}^{N}V^{i}_{l}\bigg].
\end{align},

where $\sigma_{i}$ is the $i^{th}$ eigenvalue of $M_{k}$ and $V^{i}_{l}$ is the $l^{th}$ component of the associated eigenvector, and $N=size(C_{k})$ 
\item Compute the distance of $\widehat{X}$ to each POD feature subspace using Formula 
\begin{equation}
dis(\widehat{Y}, F) = \sqrt{ \Vert \widehat{Y}\Vert ^{2} - \sum_{i=1}^{N}\alpha_{i}^{2} }.
\end{equation}
In the above,
\begin{equation}
\Vert \widehat{Y}\Vert ^{2} = \widehat{Y}\widehat{Y}= \phi(\widehat{X})\phi(\widehat{X}) = K(\widehat{X},\widehat{X}).
\end{equation}
\item \textbf{Decision}: Assign $\widehat{X}$ to the class corresponding to the minimum computed distance.
\end{enumerate}

\section{RBF parameter Formula}\label{3}

The proposed formula is based on the simultaneous minimization of class size and the maximization of inter-class distances in the feature space. We refer to this formulation as the dual min-max (DMM) formula. As a measure of class size, we consider the class diameter.

We begin by introducing the necessary definitions and notations.
\subsection{Definitions and Notations}
From a mathematical point of view, the diameter of a set $A$ in a metric space $M$ is defined as
\begin{equation*}
\mathrm{diam}(A) = \max_{x,y \in A} \|x - y\|.
\end{equation*}
The distance between two subsets $A$ and $B$ is defined as
\begin{equation*}
\mathrm{dist}(A,B) = \min_{x \in A,\, y \in B} \|x - y\|.
\end{equation*}

Denote by $D_k$ the diameter of a given class $C_k$ in the attribute space:
\begin{equation*}
D_k = \max_{X_i, X_j \in C_k} \| X_j - X_i \|.
\end{equation*}
Similarly, $\overline{D}_k$ denote the diameter of class $C_k$ in the feature space:
\begin{equation*}
\overline{D}_k = \max_{X_i, X_j \in C_k} \| \varPhi(X_j) - \varPhi(X_i) \|.
\end{equation*}

Next, denote by $d_{l,k}$ the distance between two classes $C_l$ and $C_k$ in the attribute space, and by $\overline{d}_{l,k}$ the corresponding distance in the feature space:
\begin{equation*}
d_{l,k} = \min_{X_i \in C_l,\, Y_j \in C_k} \| X_i - Y_j \|,
\end{equation*}
and
\begin{equation*}
\overline{d}_{l,k} = \min_{X_i \in C_l,\, Y_j \in C_k} \| \varPhi(X_i) - \varPhi(Y_j) \|.
\end{equation*}

\subsection{Objective and Formula Derivation}
As stated above, the objective is to select the mapping $\varPhi$ such that the class diameters 
$\overline{D}_k$ are \emph{minimized} while the inter-class distances $\overline{d}_{l,k}$ are 
\emph{maximized} in the feature space, thereby achieving optimal data separability.

First, we calculate $\overline{D}_{k}$ and $\overline{d}_{l,k}$ for the RBF kernel.  
We have:

\begin{gather*}
 \| \varPhi(X_{j}) - \varPhi(X_{j}) \|^{2} = ( \varPhi(X_{j}) - \varPhi(X_{j}), \varPhi(X_{j}) - \varPhi(X_{j})  ) = \\  \varPhi(X_{i})\varPhi(X_{i}) +  \varPhi(X_{j}) \varPhi(X_{j}) -2\varPhi(X_{i})\varPhi(X_{j}) = K(X_{i}, X_{i}) + K(X_{j}, X_{j}) - 2K(X_{i}, X_{j})
\end{gather*}

\noindent

Substituting the RBF kernel $K(X,Y) = e^{-\gamma \| X-Y \|^{2}} $ we obtain
\begin{align*}
	\| \varPhi(X_{j}) - \varPhi(X_{j}) \|^{2} = 2-2 e^{-\gamma \| X_{i}- X_{j} \|^{2}}
\end{align*}

Consequently,
\begin{gather*}
	\overline{D}_{k}^{2}(\gamma) = \max_{i,j} \{\| \varPhi(X_{j}) - \varPhi(Y_{j}) \|^{2} \}= \max_{i,j} ( 2-2 e^{-\gamma \| X_{i}- X_{j} \|^{2}} ) = \\ ( 2-2\min_{i,j}( e^{-\gamma \| X_{i}- X_{j} \|^{2}} )) =  ( 2-2 e^{-\max_{i,j}(\gamma \| X_{i}- X_{j} \|^{2}})  = 2-2 e^{-\gamma \max_{i,j}(\| X_{i}- X_{j} \|^{2}} = 
	2-2 e^{-\gamma D^{2}_{k}(\gamma)}.
\end{gather*}
\noindent
This shows that the diameter of any class in the feature space can be computed directly from the diameter of the class in the original attribute space.

Now calculate  $\overline{d}_{l,k}$ 

\begin{gather*}
	 \min_{i,j} \| \varPhi(X_{j}) - \varPhi(Y_{j}) \|^{2} = \min_{i,j}( \varPhi(X_{j}) - \varPhi(Y_{j}), \varPhi(X_{j}) - \varPhi(Y_{j})  ) =\min_{i,j} \\ \{ \varPhi(X_{i})\varPhi(X_{i}) +  \varPhi(Y_{j}) \varPhi(Y_{j}) -2\varPhi(X_{i})\varPhi(Y_{j}) \}=\min_{i,j}\{ K(X_{i}, X_{i}) + K(Y_{j}, Y_{j}) - 2K(X_{i}, Y_{j}) \}
\end{gather*}
Substituting the RBF kernel, we obtain
\begin{gather*}
	 \min_{i,j} \| \varPhi(X_{j}) - \varPhi(Y_{j}) \|^{2} = \min_{i,j} \{2-2e^{-\gamma \| X_{i}- X_{j} \|^{2}}\} = \\ (2-2\max_{i,j}e^{-\gamma \| X_{i}- X_{j} \|^{2}}) = (2-2e^{-\gamma \min_{i,j}\| X_{i}- X_{j} \|^{2}}) = (2-2e^{-\gamma d^{2}_{l,k}})
\end{gather*}
Finally
\begin{gather*}
	\overline{d}^{2}_{l,k}(\gamma) = (2-2e^{-\gamma d^{2}_{l,k}})
\end{gather*}

Similarly, the inter-class distance in the feature space can be expressed in terms of the inter-class distance in the original attribute space.

Note that, to minimize the diameters of all attribute classes, it is sufficient to minimize the largest diameter. Likewise, to maximize the inter-class distances, it is sufficient to maximize the smallest inter-class distance.

That is,

 \begin{gather*}
 \overline{D}_{max}^{2}=	\max_{k}\{\overline{D}_{k}^{2}(\gamma)\}  = 	\max_{k}\{2-2 e^{-\gamma D^{2}_{k}}\} =  2-2 \min_{k}\{e^{-\gamma D^{2}_{k}}\}=2-2 e^{-\gamma  \max_{k}\{D^{2}_{k}\}}) = \\ =2-2 e^{\gamma  D_{max}^{2}}
 \end{gather*}

And 

\begin{gather}\label{d2min}
\overline{d}_{min}^{2}(\gamma) = \min_{l,k}\{\overline{d}^{2}_{l,k}(\gamma)\} = \min_{l,k}\{2-2e^{-\gamma d^{2}_{l,k}}\} = 2-2 \max_{l,k}\{e^{-\gamma d^{2}_{l,k}}\} = 2-2e^{-\gamma \min_{l,k}\{d^{2}_{l,k}\}} = \nonumber  \\ 2-2e^{-\gamma d_{min}^{2}} 
\end{gather}

Now we need to simulatneousely minimize $\overline{D}^{2}_{max}(\gamma)$ and maximize $\overline{d}^{2}_{min}(\gamma)$ as functions of $\gamma$.

We have
\begin{equation}
\min_{\gamma}(\overline{D}^{2}_{max}(\gamma)) = 2-2e^{-\min_{\gamma}\gamma D^{2}_{max}},
\end{equation}
and

\begin{equation}
\max_{\gamma}(\overline{d}^{2}_{min}(\gamma)) = 2-2e^{-\max_{\gamma}\gamma d^{2}_{min}},
\end{equation}

This amounts to minimizing the function

\begin{equation}
F_{1}(\gamma) = \gamma D^{2}_{max} ,
\end{equation}

and maximize
\begin{equation}
F_{2}(\gamma) = \gamma d^{2}_{min}.
\end{equation}

Note that maximizing $F_{2}$ is equivalent to minimize $1/F_{2}$

Finally, we need to minimize simultaneously,

\begin{equation}\label{minmaxG}
G_{1}(\gamma) = \gamma D^{2}_{max} 
\end{equation}
and
\begin{equation}
G_{2}(\gamma) = \dfrac{1}{\gamma d^{2}_{min}} 
\end{equation} 

This can be accomplished using the weighted-form technique by defining the function

\begin{equation}\label{G(gamma)}
G(\gamma) = \lambda G_{1}(\gamma) +(1-\lambda) G_{2}(\gamma)=\lambda \gamma D^{2}_{max}  + (1-\lambda)\dfrac{1}{\gamma d^{2}_{min}} 
\end{equation}

Next, we look for the Pareto optima.  
To do so, we calculate the derivative of $G(\gamma)$ with respect to $\gamma$:

\begin{equation}
G'(\gamma) = \lambda D^{2}_{max} -(1-\lambda)\dfrac{1}{\gamma^{2}d^{2}_{min}}
\end{equation}
Then solve
\begin{equation*}
G'(\gamma) = 0
\end{equation*}
We obtain

\begin{equation}
\overline{\gamma}(\lambda)= \sqrt{\dfrac{(1-\lambda)}{\lambda D^{2}_{max}d^{2}_{min}}}
\end{equation}
 We Calculate the second derivative of $G(\gamma)$ with respect to $\gamma$:
 \begin{equation}
G"(\gamma) = \dfrac{2(1-\lambda)}{d^{2}_{min}\gamma^{3}} > 0
\end{equation}
We deduce that $G$ is a concave up function, and hence $\overline{\gamma}(\lambda)$ corresponds to a true minimums for all $\lambda \in [0,1]$.

Following this, we find the smallest Pareto minimum by substituting $\overline{\gamma}(\lambda)$ into the expression of $G(\gamma)$ \eqref{G(gamma)} to obtain all minima.  
We then minimize the resulting function with respect to $\lambda$. That is,

\begin{equation}
G(\overline{\gamma}(\lambda)) = \lambda \sqrt{\dfrac{(1-\lambda)}{\lambda D^{2}_{max}d^{2}_{min}}} D^{2}_{max}  + (1-\lambda)\dfrac{1}{\sqrt{\dfrac{(1-\lambda)}{\lambda D^{2}_{max}d^{2}_{min}}} d^{2}_{min}} =2 \sqrt{\lambda(1-\lambda)} \dfrac{D_{max}}{d_{min}}
\end{equation}

Differentiate $G$ with respect to $\gamma$

\begin{equation}
\frac{dG(\overline{\gamma}(\lambda))}{d \lambda} = \dfrac{-2\lambda +1}{\sqrt{\lambda(1-\lambda)}}\dfrac{D_{max}}{d_{min}}
\end{equation}.

Setting
\begin{equation}
\frac{dG(\overline{\gamma}(\lambda))}{d \lambda} = 0
\end{equation}
yields
\begin{equation}
\lambda = \frac{1}{2},
\end{equation}
which corresponds to a balanced weighting.

Substituting in the expression of $\overline{\gamma}(\lambda)$ we obtain the desired formula for $\gamma$.

\begin{equation}\label{gamma}
\boxed{\overline{\gamma}= \frac{1}{D_{max}d_{min}} }
\end{equation}
or
\begin{equation}\label{sigma}
\boxed{\overline{\sigma} = \sqrt{\frac{1}{2}D_{max}d_{min}} }
\end{equation}

In practice, if an attribute vector is mistakenly located inside the wrong class and lies very close to the vectors of that class, $d_{\min}$ may be calculated incorrectly and approach zero. To avoid such situations, it is preferable to relax the formula and use the average of the inter-class distances instead of the minimum. Therefore, $d_{\min}$ is replaced by the average distance $d_{\rm av}$, defined as

\begin{equation*}
d_{av}= \sqrt{ \frac{1}{T}\sum_{\substack{i,j  \\ i <j }}d^{2}_{i,j} },
\end{equation*}
where $T$ denotes the total number of distances $d_{i,j}^{2}$.

Consequently
\begin{equation}\label{gamma1}
\boxed{\overline{\gamma}= \frac{1}{D_{max}d_{av}} }
\end{equation}
and
\begin{equation}\label{sigma1}
\boxed{\overline{\sigma} = \sqrt{\frac{1}{2}D_{max}d_{av}} }
\end{equation}
Tests show that, for a small number of classes, both formulas achieve the same accuracy. However, for a large number of classes, the formula using the average inter-class distances provides slightly better results. Therefore, in the validation section, the latter formula will be used.

\section{Validation Tests}
\label{4} 
The proposed DMM formula is verified and validated on selected real-world datasets from the LIBSVM~\cite{LIBSVM} and OpenML~\cite{openML} repositories, using both SVM and KOS classifiers. Results are compared against those obtained with the standard learning process for both methods. Note that for SVM, the penalty parameter $C$ still needs to be tuned, while the parameter $\gamma$ is determined by the DMM formula.

The LIBSVM software package is used for SVM, whereas KOS is implemented according to the process described in \cite{lremaki2026}.

Tables~\ref{Leaen-DMM comp SVM} and ~\ref{Leaen-DMM comp KOS} show the accuracy and precision for both SVM and KOS using the standard learning process and the proposed DMM formula. The results indicate that the DMM formula achieves comparable performance to the standard learning process, demonstrating the effectiveness of the proposed formula. Table~\ref{Processing-time} presents the processing times for the tests, showing, as expected, a significant reduction in computational time.

\begin{table}[H]
\caption{SVM: Learning Performance vs. DMM formula Comparison (Accuracy|Precision).\label{Leaen-DMM comp SVM}}
	\begin{adjustwidth}{-\extralength}{0cm}
		\begin{tabularx}{\fulllength}{p{3cm}C}
			\toprule
			\textbf{Name/Classes}& \textbf{SVM \linebreak  Learning ~~--~~ DMM \linebreak \mbox{\small (Acc|Prec)~~--~~(Acc|Prec)}}  \\
			\midrule
{Leukemia/2}	
 & $(82.35\%|86.42\%) -(\textbf{82.35\%|86.42\%}) $ \\                    
                   \midrule
{svmguide1/2}	 & $(96.87\%|96.87\%)-(\textbf{96.87\%|96.88\%})$ \\        
                   \midrule
{splice/2}	 & $(90.43\%|90.50\%)-(\textbf{89.66\%|89.71\%})$  \\          
                   \midrule
{austrian/2} & $(85.04\%|85.01\%)-(\textbf{81.52\%|81.50\%})$  \\          
                   \midrule
{madelon/2}	 & $(61.16\%|59.51\%)-(\textbf{60.16\%|60.25\%})$ \\         
                   \midrule
{DNA/3} & $(94.43\%|94.14\%)-(\textbf{95.62\%|95.64\%})$  \\             
                   \midrule
{Satimage/6}	 & $(91.85\%|91.68\%)-(\textbf{90.40\%|90.23\%})$ \\           
                   \midrule
{USPS/10}	& $(95.26\%|95.32\%)-(\textbf{95.72\%|95.75\%})$ \\       
                   \midrule
{letter/26}	 & $(97.9\%|97.90\%)-(\textbf{97.32\%|97.33\%})$ \\    
                   \midrule
{shuttle/7}	 & $(99.93\%|99.93\%)-(\textbf{99.88\%|99.88\%})$ \\
  				   \midrule
{BreastCancer/2}	 & $(97.07\%|97.07\%)-(\textbf{97.66\%|97.60\%})$  \\
  				   \midrule 
{Zernike/10}	 & $(82.00\%|82.03\%)-(\textbf{82.33\%|82.47\%})$ \\
  				   \midrule 	
{Diabetes/2}	 & $(74.026\%|74.40\%)-(\textbf{71.42\%|70.51\%})$   \\
  				   \midrule 	
{mfeat-morphological/10}	 & $(73.66\%|69.36\%)-(\textbf{73.33\%|74.70\%})$  \\
           \bottomrule
\end{tabularx}
\end{adjustwidth}
\end{table}

\begin{table}[H]
\caption{KOS: Learning Performance vs. DMM formula Comparison (Accuracy|Precision).\label{Leaen-DMM comp KOS}}
	\begin{adjustwidth}{-\extralength}{0cm}
		\begin{tabularx}{\fulllength}{p{3cm} C}
			\toprule
			\textbf{Name/Classes} & \textbf{KOS \linebreak  Learning ~~--~~ DMM \linebreak \mbox{(Acc|Prec)~~--~~(Acc|Prec)}} \\
			\midrule
{Leukemia/2}	& $(82.35\%|82.58\%) - (\textbf{82.35\%|82.44\% })$ \\                    
                   \midrule
{svmguide1/2}	& $(96.30\%|96.31\%)-(\textbf{96.37\%|96.38\%})$ \\        
                   \midrule
{splice/2}	 & $(89.42\%|89.41\%)-(\textbf{89.33\%|89.31\%})$ \\          
                   \midrule
{austrian/2} & $(85.63\%|85.36\%)-(\textbf{82.11\%|81.80\%})$ \\          
                   \midrule
{madelon/2}	 & $(61.50\%|61.55\%)-(\textbf{56.16\%|56.21\%})$ \\         
                   \midrule
{DNA/3}  & $(92.32\%|92.83\%)-(\textbf{91.48\%|93.74\%})$ \\             
                   \midrule
{Satimage/6}	 & $(91.35\%|90.40\%)-\textbf{90.40\%|90.40\%})$ \\           
                   \midrule
{USPS/10}	& $(95.76\%|95.76\%)-(\textbf{95.86\%|95.59\%})$ \\       
                   \midrule
{letter/26}	 & $(97.42\%|97.44\%)-(\textbf{96.38\%|96.45\%})$ \\    
                   \midrule
{shuttle/7}	 & $(99.91\%|99.91\%)-(\textbf{99.66\%|99.66\%})$ \\
  				   \midrule
{BreastCancer/2} & $(97.07\%|97.01\%)-(\textbf{97.66\%|97.82\%})$ \\
  				   \midrule 
{Zernike/10}	 & $(81.00\%|80.82\%)-(\textbf{80.16\%|80.20\%})$ \\
  				   \midrule 	
{Diabetes/2}	  & $(74.46\%|72.20\%)-(\textbf{70.56\%|70.56\%})$ \\
  				   \midrule 	
{mfeat-morphological/10} & $(73.833\%|75.07\%)-(\textbf{72.33\%|72.61\%})$ \\
           \bottomrule
\end{tabularx}
\end{adjustwidth}
\end{table}


\begin{table}[H]
\caption{DMM formula vs. Learning processing-time~comparison.\label{Processing-time}}
	\begin{adjustwidth}{-\extralength}{0cm}
		\begin{tabularx}{\fulllength}{p{3cm}CC}
			\toprule
			\textbf{Name/ Classes} & \textbf{SVM \linebreak Learning ~~--~~ DMM} & \textbf{KOS \linebreak Learning ~~--~~ DMM} \\
			\midrule
{Leukemia/2}	 & $0m8.646s - \textbf{0m0.0382s}$& $0m0.005s - \textbf{0m0.0041s}$ \\
                   \midrule
{svmguide1/2} & $0m33.327s - \textbf{0m2.280s}$ & $0m1.681s - \textbf{ 0m0.634
s }$ \\
                   \midrule
{splice/2}	 & $0m29.190s - \textbf{0m2.720s}$ & $0m0.263s - \textbf{0m0.0761s }$ \\
                   \midrule
{austrian/2}  & $0m1.721s -\textbf{0m0.228s }$ & $0m0.0260s - \textbf{0m0.006s} $ \\
                   \midrule
{madelon/2} & $16m50.721s -\textbf{1m33.509s }$ & $0m1.571s - \textbf{0m0.243s} $ \\
                   \midrule
{DNA/3} & $4m24.300s -\textbf{0m16.491s }$ & $0m0.614s - \textbf{0m0.115s} $ \\
                   \midrule
{Satimage/6}	 & $3m53.042s -\textbf{0m11.003s }$& $0m1.468s - \textbf{0m0.242s} $ \\
                   \midrule
{USPS/10}	  & $81m18.931s - \textbf{2m39.773s}$ & $0m6.347s - \textbf{0m0.754s} $ \\
                   \midrule
{letter/26} & $44m11.996s - \textbf{2m5.038s }$ & $\textbf{0m21.330s - 0m1.904s} $ \\
                   \midrule
{shuttle/7}	 & $62m28.783s - \textbf{1m58.975}$ & $0m55.077s - \textbf{0m6.120s} $ \\
 \midrule
{BreastCancer/2}	 & $0m1.696s - \textbf{0m0.630s}$ & $0m0.037s - \textbf{0m0.006s}$ \\                 
                   \midrule
{Zernike/10}	 & $0m48.875s - \textbf{0m3.718s }$ & $0m0.127s - \textbf{0m0.0190s}$ \\        
                   \midrule
{Diabetes/2}	 & $0m6.959s - \textbf{0m0.729s}$ & $0m0.062s - \textbf{0m0.007s}$ \\          
                   \midrule
{mfeat-morphological/10}  & $0m23.421s - \textbf{0m4.200s }$ & $0m0.1069s - \textbf{0m0.0160s}$ \\          
                   \midrule
		\end{tabularx}
	\end{adjustwidth}
\end{table}

\vspace{-9pt}

\section{Conclusions}\label{5}
In this paper, we propose an analytical formula to estimate the RBF kernel parameter $\gamma$, which is extensively used in kernel methods for data classification. This parameter must be carefully tuned, as the performance of these methods strongly depends on its value and is highly sensitive to improper selection. In practice, $\gamma$ is typically determined through computationally expensive training procedures, which may also lack robustness.  

The proposed formula is derived by optimizing the diameters of mapped classes in the feature space while simultaneously maximizing inter-class feature distances. The mathematical details are presented, and tests are performed using both SVM and KOS methods, comparing the standard learning process with the proposed formula. The results show comparable performance, while the computational time is significantly reduced.  

Moreover, using an analytical formula has the advantage of leveraging the entire training set, whereas standard tuning techniques require splitting the data into training and validation sets, which can negatively affect the robustness of the process.

\begin{adjustwidth}{-\extralength}{0cm}
\end{adjustwidth}

\begin{thebibliography}{999}



\bibitem{Ethem20} 
Alpaydin, E. \emph{Introduction to Machine Learning}, 4th ed.; MIT:  Cambridge, MA, USA, 
 2020; pp. xix, 1–3, 13–18
.

\bibitem{Vap95} 
Vapnik, V. \emph{The Nature of Statistical Learning Theory}; Springer: New York, NY, USA, 1995.

\bibitem{Vap98} 
Vapnik, V. \emph{Statistical Learning Theory}; Wiley: New York, NY, USA, 1998.

\bibitem{CorVap95} 
Cortes, C.; Vapnik, V. Support vector networks. \textit{Mach. Learn.} \textbf{1995}, \emph{20}, 273--297.

\bibitem{Merc909} Mercer, J. Functions of positive and negative type and their
connection with the theory of integral equations. \textit{Philos. Trans. R. Soc.} \textbf{1909}, \emph{209}, 415–446.

\bibitem{Wang023} Wang, H.; Li, G.; Wang, Z. Fast SVM classifier for large-scale classification problems. \textit{Inf. Sci.} \textbf{2023}, \emph{642}, 119136.

\bibitem{Shao023} Shao, Y.H.; Lv, X.J.; Huang, L.W.; Bai, L. Twin SVM for conditional probability estimation in binary and multiclass classification. \textit{Pattern Recognit.} \textbf{2023}, \emph{136}, 109253.

\bibitem{Wang024} Wang, H.; Shao, Y. Fast generalized ramp loss support vector machine for pattern classification. \textit{Pattern Recognit.} \textbf{2024}, \emph{146}, 109987.

\bibitem{BWang2021} Wang, B.Q.; Guan, X.P.; Zhu, J.W.; Gu, C.C.; Wu, K.J.; Xu, J.J. SVMs multi-class loss feedback based discriminative dictionary learning for image classification. \textit{Pattern Recognit.} \textbf{2021}, \emph{112}, 107690.

\bibitem{Borah020} Borah P.; Gupta, D. Functional iterative approaches for solving support vector classification problems based on generalized Huber loss. \textit{Neural Comput. Appl.} \textbf{2020}, \emph{32}, 1135–1139.

\bibitem{Gaye021} Gaye, B.; Zhang, D.; Wulamu, A. Improvement of Support Vector Machine Algorithm in Big Data Background. \textit{Hindawi Math. Probl. Eng.} \textbf{2021}, \emph{2021}, 5594899. https://doi.org/10.1155/2021/5594899.


\bibitem{Tian012} Tian, Y.; Shi, Y.; Liu, X. Advances on support vector machines research. \textit{Technol. Econ. Dev. Econ.} \textbf{2012}, \emph{18}, 5–33. https://doi.org/10.3846/20294913.2012.661205.

\bibitem{KMOD001} Ayat, N.E.; Cheriet, M.; Remaki, L.; Suen, C.Y. KMOD---A New Support Vector Machine Kernel with Moderate Decreasing for Pattern Recognition. Application to Digit Image Recognition. In Proceedings of the Sixth International Conference on Document Analysis and Recognition, Seattle, WA, USA, 10--13 September 2001; pp. 1215--1219.

\bibitem{lremaki2026} Remaki, L. Kernel-Based Optimal Subspaces (KOS): A Method for Data Classification. In \textit{Machine Learning and Knowledge Extraction}; Basel Vol. 8, Iss. 2,  (2026): 52. DOI:10.3390/make8020052 .

\bibitem{lremaki2024} Remaki, L. Efficient Alternative to SVM Method in Machine Learning. In \textit{Intelligent Computing}; Arai, K., Ed.; Lecture Notes in Networks and Systems; Springer: Cham, Switzerland, 2025; Volume 1426.

\bibitem{subspace2016} Yoshikazu, W.; Nakayama, Y. Learning subspace classification using subset approximated kernel principal component analysis. \textit{IEICE Trans. Inf. Syst.} \textbf{2016}, \emph{99}, 1353--1363.

\bibitem{subspace2022} Jiang, W.; Chen, Y.; Wu, L.; Yu, P.S. Subspace learning for effective meta-learning. In Proceedings of the 39th International Conference on Machine Learning, Baltimore, MD, USA, 17--23 July 2022; pp. 10177--10194.

\bibitem{subspace2023} Cao, Y.-H.; Wu, J.X. Random subspace sampling for classification with missing data. \textit{J. Comput. Sci. Technol.} \textbf{2024}, \emph{39}, 472-486.


\bibitem{kpca1997} Schölkopf, B.; Smola, A.; Müller, K.R. Kernel principal component analysis. In \textit{Artificial Neural Networks—ICANN’97}; Springer: Berlin, Germany, 1997; pp. 583–588.

\bibitem{ksubspace2022} Zhou, S.; Ou, Q.; Liu, X.; Wang, S.; Liu, L.; Wang, S. Multiple Kernel Clustering with Compressed Subspace Alignment. \textit{IEEE Trans. Neural Netw. Learn. Syst.} \textbf{2021}, \emph{33}, 252-263.

\bibitem{LIBSVM} Chang, C.C.; Lin, C.J. LIBSVM: A library for support vector machines. \textit{ACM Trans. Intell. Syst. Technol.} \textbf{2011}, \emph{2}, 27:1--27:27
.

\bibitem{Lin2008} Lin, S.W.; Ying, K.C.; Chen, S.C.; Lee, Z.J. Particle swarm optimization for parameter determination and feature selection of support vector machines. \textit{Expert Syst. Appl.} \textbf{2008}, \emph{35}, 1817--1824.

\bibitem{Syarif2016} Syarif, I.; Prugel-Bennett, A.; Wills, G. SVM Parameter Optimization Using Grid Search and Genetic Algorithm to Improve Classification Performance. \textit{TELKOMNIKA} \textbf{2016}, \emph{14}, 1502--1509
. https://doi.org/10.12928/TELKOMNIKA.v14i
.

\bibitem{Shekar19} Shekar, B.H.; Dagnew, G. Grid Search-Based Hyperparameter Tuning and Classification of Microarray Cancer Data. In Proceedings of the Second International Conference on Advanced Computational and Communication Paradigms (ICACCP), Gangtok, India, 25--28 February 2019.

\bibitem{Hinton06} Hinton, E.; Osindero, S.; Teh, Y. A fast learning algorithm for deep belief nets. \textit{Neural Comput.} \textbf{2006}, \emph{18}, 1527–1554.

\bibitem{Bergstra12} Bergstra, J.; Bengio, Y. Random Search for Hyper-Parameter Optimization. \textit{J. Mach. Learn. Res.} \textbf{2012}, \emph{13}, 281--305.

\bibitem{volkwein2013proper} Volkwein, S. \emph{Proper Orthogonal Decomposition: Theory and Reduced-Order Modelling}; Lecture Notes; University of Konstanz: Konstanz, Germany, 2013; Volume 4.

\bibitem{Wang017} Wang, W.; Zhang, M.; Wang, D.; Jiang, Y. Kernel PCA feature extraction and the SVM classification algorithm for multiple-status, through-wall, human being detection. \textit{EURASIP J. Wirel. Commun. Netw.} \textbf{2017}, \emph{2017}, 151. https://doi.org/10.1186/s13638-017-0931-2.

\bibitem{SMOTE} Chawla, N.V.; Bowyer, K.W.; Hall, L.O.; Kegelmeyer, W.P. Smote: Synthetic minority over-sampling technique. \textit{J. Artif. Intell. Res.} \textbf{2002}, \emph{16}, 321–357.

\bibitem{Remaki2000} Remaki, L.; Cheriet, M. KCS---New kernel family with compact support in scale space: Formulation and impact. \textit{IEEE Trans. Image Process.} \textbf{2000}, \emph{9}, 970--981.

\bibitem{Koenderink1984} Koenderink, J.J. The structure of images. \textit{Biol. Cybern.} \textbf{1984}, \emph{53}, 363–370.

\bibitem{openML} Vanschoren, J.; Van Rijn, J.N.; Bischl, B.; Torgo, L. 
OpenML: Networked science in machine learning.
\textit{SIGKDD Explor. Newsl.} \textbf{2014}, \emph{15}, 49–60.

\bibitem{Leu99} Golub, T.R.; Slonim, D.K.; Tamayo, P.; Huard, C.; Gaasenbeek, M.; Mesirov, J.P.; Coller, H.; Loh, M.L.; Downing, J.R.; Caligiuri, M.A.; et al. Molecular classification of cancer: Class discovery and class prediction by gene expression monitoring. \textit{Science} \textbf{1999}, \emph{286}, 531.

\bibitem{svmguide103} Hsu, C.W.; Chang, C.C.; Lin, C.J. \textit{A Practical Guide to Support Vector Classification}; Technical Report; Department of Computer Science, National Taiwan University: Taipei, Taiwan, 2003.

\bibitem{splice} Available online: http://archive.ics.uci.edu/ml/index.php 

\bibitem{austrian} Available online: https://www.csie.ntu.edu.tw/cjlin/libsvmtools/datasets/.

\bibitem{madelon05} Guyon, I.; Gunn, S.; Ben Hur, A.; Dror, G. Result analysis of the NIPS 2003 feature selection challenge. \textit{Adv. Neural Inf. Process. Syst.} \textbf{2004, pp. 545–552}, \emph{17}
.

\bibitem{DNA02} Hsu, C.W.; Lin, C.J. A comparison of methods for multi-class support vector machines. \textit{IEEE Trans. Neural Netw.} \textbf{2002}, \emph{13}, 415–425.

\bibitem{USPS94} Hull, J.J. A database for handwritten text recognition research. \textit{IEEE Trans. Pattern Anal. Mach. Intell.} \textbf{1994}, \emph{16}, 550–554.


\bibitem{tensorflow2016} Abadi, M.; Barham, P.; Chen, J.; Chen, Z.; Davis, A.; Dean, J.; Devin, M.; Ghemawat, S.; Irving, G.; Isard, M.;~et~al. TensorFlow: A System for Large-Scale Machine Learning. In Proceedings of the 12th USENIX Symposium on Operating Systems Design and Implementation, OSDI 16, Savannah, GA, USA, 2--4 November 2016; pp. 265--283.

\bibitem{MatlabCplusplus} Andrews, T. Computation Time Comparison Between Matlab and C++ Using Launch Windows. Available online: \text{https://digitalcommons.calpoly.edu/aerosp/78/?utm\_source=chatgpt.com.}

\bibitem{Kernels2014} Kung, S. Y. Kernel Methods and Machine Learning. Cambridge University Press, 2014.
\bibitem{Kernels2004} Shawe-Taylor, J., Cristianini, N. Kernel Methods for Pattern Analysis. Cambridge University Press, 2004.
\end{thebibliography}
\end{document}